\documentclass[11pt]{article}

\usepackage{amssymb}
\usepackage{times}
\usepackage{url}
\usepackage{hyperref}

\newcommand{\F}{\mathbb{F}}
\newcommand{\Q}{\mathbb{Q}}
\newcommand{\C}{\mathbb{C}}
\newcommand{\GL}{\mathrm{GL}}
\newcommand{\SL}{\mathrm{SL}}
\newcommand{\PSL}{\mathrm{PSL}}
\newcommand{\GF}{\mathrm{GF}}
\newcommand{\glnf}{\mathrm{GL}(n,\F)}

\newcommand{\abk}{\allowbreak}

\begin{document}

\title{L.~G.~Kov\'{a}cs and linear groups}

\author{A.~S. Detinko and D.~L. Flannery}

\date{}

\maketitle

\bigskip

\begin{center}

Dedicated to the memory of Laci Kov\'{a}cs

\end{center}

\bigskip

\begin{abstract}
We survey the legacy of L.~G.~Kov\'{a}cs in linear group
theory, with a particular focus on classification questions.
\end{abstract}

\bigskip

\medskip

Classifying linear groups is an old problem.
For given degree $n>1$, field $\F$, and group type, the task is to list 
irredundantly all subgroups of $\glnf$ of that type up to $\glnf$-conjugacy. 
If possible, each conjugacy class representative should be specified
by a generating set of matrices. Examples of the group type are: 
finite, soluble, nilpotent, quasi-simple, 
maximal in its class (provided that each relevant subgroup 
of $\glnf$ is contained in a maximal), irreducible, primitive, 
monomial, generated by matrices with special properties
(e.g., pseudoreflections), and so on.
 
To narrow the scope and thereby have a reasonable hope of solving
a linear group classification problem, we impose extra conditions,
such as the following.
\begin{itemize}
\item Characteristic: zero ($\F = \C$, the complex field, or $\F =
\Q$, the rationals, are typical instances); positive (mostly
finite $\F$).
\item  Degree: `small', or otherwise restricted, 
according to the prime factorisation of $n$.
\end{itemize}
Furthermore, although $\GL(n,\F)$-conjugacy is a natural
classification criterion, we might ask for isomorphism class 
representatives instead. We may even limit ourselves to 
classifying subgroups of $\SL(n,\F)$ or $\mathrm{PSL}(n,\F)$.

Laci Kov\'{a}cs had an abiding interest in linear group theory.
He was one of the first to realise the suitability of computer algebra
systems as an environment for linear group classification: using 
a computer to aid in the compilation of lists, and applying implemented
lists to solve related algorithmic problems. 

Much of Laci's research dealing with representation theory and
permutation groups has strong intersections with linear group
theory. Especially pertinent here are the asymptotic bounds that
he proved for finite soluble and nilpotent groups. We 
mention a few of these results, to give context; the paper by
G.~R.~Robinson in this volume contains more detail.

In \cite{KoLG1}, Laci and J.~D.~Dixon derived a bound in terms of
$n$ and $\F$ on the number of generators of a finite nilpotent subgroup
 of $\glnf$, where $\F$ is a finite degree extension of its prime 
subfield. Then Laci, R.~M.~Bryant, and Robinson extended that 
result to any finite group generated by its soluble radical and 
generalised Fitting subgroup~\cite{KoLG2}. The case of soluble linear 
groups was crucial (note also the paper~\cite{KoLGSim} with
H.-S.~Sim on the number of generators of an abstract
finite soluble group). But perhaps the most striking
achievement in this area is \cite{KoLG4}. Laci and Robinson prove 
that a finite completely reducible linear group of degree $n$ over 
any field can be generated by $\left\lfloor \frac{3n}{2}\right\rfloor$ 
elements. Further contributions to linear group theory include
\cite{KoLG6}, which establishes complete reducibility of
representations of the monoid of $n\times n$ matrices over a
finite field.

Below we discuss a research programme founded by Laci in the 
1990s, that targeted difficult linear group classification
problems. We describe how the objectives of this
programme were carried out with some of his students and postdoctoral
researchers. Laci guided the development of techniques and 
formulated major strategies in this programme.

\section{Background}

We begin with a sketch of historical background.
See \cite[Chapter~III]{ZalesskiiI} and \cite[Chapter~3,
\S4]{ZalesskiiII} for comprehensive surveys.

Early interest in soluble linear groups over finite fields stemmed
from their connection to soluble permutation groups. C.~Jordan
gave a method (which can be viewed as an archetypal group-theoretic
algorithm) to construct such linear groups from those of smaller
degree. Jordan's treatment is cumbersome and does not give
a full classification up to conjugacy.
For finite groups in characteristic zero, it essentially
suffices to classify the irreducible ones, dispensing 
with small degrees first. Degrees $2$ and $3$ were investigated 
by Jordan and F.~Klein, amongst others.
H.~F.~Blichfeldt's wonderful book~\cite{Blichfeldt} covers finite 
complex linear groups of degree at most $4$. 

It became customary to ignore imprimitive groups. 
The emphasis was rather on primitive subgroups of $\SL(n,\C)$ or 
their images in $\PSL(n,\C)$, and sometimes groups were  
determined only up to isomorphism. The standard justification
for this is as follows. Let $\F$ be
algebraically closed, and $Z$ be the scalar subgroup of
$\GL(n,\F)$. Given $G\leq \GL(n,\F)$ we may define $H\leq \SL(n,\F)$  
such that $GZ=\abk HZ$. Suppose that $G$ is irreducible (resp., primitive).
Then $H$ is irreducible (resp., primitive), $H$ and $G$ have isomorphic 
central quotients, and $H$ is finite precisely when $G/Z(G)$ is finite. 
An advantage of this reduction is that there are only finitely many 
conjugacy classes of finite primitive subgroups of $\SL(n,\F)$. 
However, as W.~Feit has pointed out~\cite{Feitpersonal}, producing a 
classification in $\GL(n,\C)$ from one in $\SL(n,\C)$ or $\PSL(n,\C)$ 
is not at all straightforward. \footnote{`Don't
fall into the trap I fell into, Blichfeldt only classified the
groups in dimension $4$ in $PGL(4,\bf{C})$. It is a long way to go
to $GL(4,\bf{C})$.'}

Classifying finite primitive (or quasiprimitive) subgroups of
$\SL(n,\C)$ gained popularity in the lead-up to the classification
of finite simple groups.
After Blichfeldt, authors including R.~Brauer, Feit, W.~C.~Huffman,
J.~H.~Lindsey II, and D.~B.~Wales gave accounts for $n \leq 10$
(see \cite[pp.~76--78]{FeitSurvey} and \cite{Feit}).

Another wave of activity began in the late 1940s, as soluble
linear groups were recognised to play a fundamental role
in the theory of infinite soluble groups. 
D.~A.~Suprunenko and his students obtained various
classifications of soluble linear and permutation groups.
Usually just the maximal soluble subgroups of $\glnf$ are 
described (each soluble subgroup lies in a maximal).
For example, \cite[Theorem~6, p.~167]{Supr} classifies 
the maximal irreducible soluble subgroups of $\GL(p,q)$
up to conjugacy, $p$ prime, with an explicit generating set 
stated for each conjugacy class representative 
(cf.~\cite{DetinkoClassical}). Many of the classification results 
for soluble matrix and permutation groups by Suprunenko and his 
school are summed up in \cite{SuprPermGp}.
That book also contains a classification of the maximal
primitive soluble subgroups of $\mathrm{Sym}(n)$, 
where $n\in \{ p^q, p^{q^2}, p^{qr} \mid p, q, r \mbox{ prime} \}$.

Other notable classifications are in prime or prime-square degree: 
minimal irreducible groups~\cite[p.~2986]{ZalesskiiIII},
and irreducible $p$-groups over an algebraically closed field 
(S.~B.~Conlon). Aside from this work, and in contrast to the 
insoluble case, exhaustive classifications of soluble linear groups
were rarely attempted---until the advent of the 
research school directed by Laci.

\section{The Kov\'{a}cs School}

Some time ago, Laci, J.~Neub\"{u}ser, and M.~F.~Newman proposed
an algorithm to construct maximal subgroups of low index in a finitely
presented group~\cite[pp.~2--4]{Short}. 
Their algorithm relies on having a list of the
primitive subgroups $H$ of $\mathrm{Sym}(n)$ where $n$ is the
subgroup index. 
If $H$ is soluble then $n=\abk p^m$ for some $m$ and prime $p$, and 
listing these permutation groups is equivalent to classifying 
the irreducible soluble subgroups of $\GL(m,p)$ up to conjugacy.
The need for such information motivated the
PhD project of Laci's student Mark Short~\cite{Short}.

Short's overall approach is based on theory of maximal irreducible soluble 
subgroups of $\glnf$, as in \cite[Chapter~V]{Supr} and with 
antecedents in work of Jordan. Chapters~3--5 of \cite{Short} 
furnish a classification of the irreducible soluble subgroups of 
$\GL(2,q)$ for odd $q$ 
(A.~Hulpke later found that two conjugacy 
classes of monomial groups were missing). Other necessary results 
for $\GL(r,q)$, $r$ an odd prime, and for primitive soluble subgroups
of $\GL(4,q)$, are provided. The listing in \cite[Chapter~7]{Short} 
of imprimitive groups of degree $4$ is supplemented by a CAYLEY computation. 
Using the methods in his thesis, Short classified
the irreducible soluble subgroups of $\GL(n,p)$ up to conjugacy
for all $p^n<256$. He implemented this classification 
as a data library and made it publicly available.

A subgroup of $\GL(n,q)$ of order coprime to $q$ can be
`lifted' to an isomorphic copy in $\GL(n,\C)$, and the 
lifting respects absolute irreducibility.
Despite this link, classification problems over $\C$ have 
a different flavour to those over finite fields.
One complicating factor is that a classification of finite 
subgroups of $\GL(n,\C)$ might entail an infinite list (whereas 
there are only finitely many finite primitive subgroups of $\SL(n,\C)$ 
up to conjugacy). For the lists discussed here, we may
introduce a parametrisation on certain families of matrices, so 
that each listed group is designated by an integer string
that corresponds to a generating set made
up from the parametrised families. A model for this sort of
listing scheme is exhibited in \cite{Conlon}. Conlon 
classifies the finite irreducible $p$-subgroups of $\GL(p,\F)$,
where $\F$ is a field not of characteristic
$p$ with all $p$-power roots of unity.
Such a group is conjugate to a subgroup $G$ of the full monomial 
matrix group $C_{p^\infty}\wr C_p$.
The subgroup of diagonal matrices in $G$ has index $p$ and 
(to guarantee irreducibility) must be non-scalar.
Conlon gave presentations for groups in his list, and proved 
that any two of them are conjugate if they are isomorphic.

Attacks on other monomial group classification problems 
have followed the same basic pattern as in \cite{Conlon}. 
Suppose that $G\leq\GL(n,\C)$ is monomial. Let
$\pi$ be the natural projection of $G$ into $\mathrm{Sym}(n)$ whose
kernel $D$ is the subgroup of diagonal matrices (that is, $\pi$ sends
non-zero matrix entries to $1$).
First, all candidates for $\pi (G)$ are written down; namely 
the transitive $T\leq \mathrm{Sym}(n)$. 
Then we solve (S), the $T$-submodule problem:
find all $D$ normalised by $T$ such that if $\ker \pi = D$ 
then $G$ is irreducible. The next step is the extension problem 
(E) for each $T$ and its accompanying $T$-modules. Lastly, we 
solve the conjugacy problem (C), showing that each 
$\GL(n,\C)$-conjugacy class is represented exactly once 
in our final list.

After the success of \cite{Conlon}, composite degrees present a
new challenge. The second author, another student of Laci's, 
classified the finite irreducible linear $p$-groups of degree $p^2$ 
for $p=2$~\cite{FlanneryMemoir}. We say a little bit about the methods 
used (some of which appear in \cite{Bacskai,Short} too).
The submodule lattice of a direct sum of modules 
may be assembled from isomorphisms between sections of the summands,
via a well-known theorem due to Goursat and Remak. This is applied 
to solve (S). Second cohomology features in the solution of (E)
and (C): $|H^2(T,D)|$ is an upper bound on the number of conjugacy 
classes of $G\leq \GL(4,\C)$ such that $\pi (G) = T$ and
$\ker \pi = D$.
Lyndon-Hochschild-Serre spectral sequences are
used to calculate the requisite orders. 
For each $T$ and $D$, precisely $|H^2(T,D)|$ extensions $G$ of $D$
by $T$ in $\GL(4,\C)$ are constructed. Any remaining conjugacy 
between these extensions is eliminated by ad hoc means.
Having dealt with the $2$-groups, Flannery went on to classify
all finite irreducible monomial subgroups of
$\GL(4,\C)$~\cite{Flannery2}. 

At this juncture it is appropriate to note a question in the province 
of the submodule listing problem (S), that arose out of an algorithm
suggested by Conlon for decomposing group characters.
Let $p$ be a prime, $V$ be the central quotient of
$C_{p^\infty}\wr (C_p)^n$, and $N$ be any finite normal subgroup
of $V$. Conlon conjectured that the centre $Z$ of $V/N$ has order
at least $p^n$. Examples are known where $Z\cong\abk C_p^n$, $Z\cong
\abk C_{p^n}$, and $|Z|$ is much greater than $p^n$. When $n=1$ there
is nothing to prove. Conlon verified his conjecture for $2\leq n
\leq 4$ by a combination of hand and machine calculations.

As far as we know, Conlon's conjecture is unresolved.
Laci's student Zolt\'{a}n B\'{a}cskai made
progress towards an affirmative proof, but eventually 
changed his thesis topic to classifying finite 
irreducible monomial subgroups $G$ of $\GL(p,\C)$.
Insoluble groups now crop up: $\pi(G)=T\leq \mathrm{Sym}(p)$ is 
either `compulsory' (i.e., soluble, $\mathrm{Alt}(p)$, 
or $\mathrm{Sym}(p)$); or `sporadic', with just $11$ values of $p$ 
less than $1000$ giving such a transitive group. B\'{a}cskai obtained 
a complete classification for $p\leq 11$, and for arbitrary prime 
degree $p$ when $T$ is compulsory. In particular, the solution of (S) 
for all $T$ occupies Chapters 3 and 4 and Section~7.1 of \cite{Bacskai}
(observe that a finite soluble monomial subgroup of $\GL(p,\C)$ 
is irreducible if and only if its diagonal matrix subgroup is 
non-scalar). B\'{a}cskai's thesis, which should be in the literature, 
contains many valuable results on linear group classification. 

An irreducible linear group of prime degree is either primitive or
monomial. Dixon and A.~E.~Zalesskii classified finite primitive
subgroups of $\SL(p,\C)$,
and insoluble finite monomial groups of prime degree 
over an algebraically closed field, in~\cite{DZI,DZII,DZICorr}.
The paper \cite{DZI} has a traditional aim---classifying
primitive unimodular linear groups over $\C$---and makes critical
use of the classification of finite simple groups.

A further milestone in the Kov\'{a}cs programme was supplied by
Burkhard H\"{o}f-ling, who worked as a postdoctoral researcher with
Laci. In a long and interesting \mbox{paper}~\cite{Hoefling}, H\"{o}fling
settles the case of imprimitive non-monomial finite irreducible groups 
over $\C$ in smallest degree. 
He begins by considering the general situation of
$G\leq \abk \GL(2n,\C)$ with an unrefinable imprimitivity system of 
size $2$. Either $G$ has just one system, or it has three
(\cite[Theorem~2.4]{Hoefling} is a broader statement; its 
proof cites \cite{Ko3}). This yields an initial split in the  
classification. To construct all $G$, one needs to know all 
primitive groups of degree $n$. Thus, as part of his solution in 
degree $4$, H\"{o}fling classified the finite primitive
subgroups of $\GL(2,\C)$ up to conjugacy. 
Each such group is contained in a central product of scalars with 
a primitive subgroup of $\SL(2,\C)$, and 
the latter were classified previously~\cite[Chapter~III]{Blichfeldt}.
H\"{o}fling's lists, together with \cite{Bacskai,Flannery2} 
and \cite[Chapters~V, VII]{Blichfeldt} filled out to $\GL(n, \C)$,
would complete the Kov\'{a}cs programme in degrees less than 
$5$ over $\C$. Degree $5$ is surely achievable too, with 
\cite{Bacskai,DZI} as a foundation.

Laci's student Hyo-Seob Sim wrote several papers on
metacyclic linear groups. In \cite{Sim}, he
examines the structure of metacyclic primitive subgroups of
$\GL(n,q)$, with the intent to classify these groups when 
$n$ is an odd prime power (cf.~the classification \cite{Glasgow} 
of nilpotent primitive subgroups of $\GL(n,q)$ for all $n$, $q$; 
and the {\sf GAP} procedure in \cite{Nilmat} that returns the 
groups for input $n$, $q$). Laci and Sim in \cite{KoLG5} give a 
condition to decide whether two nilpotent metacyclic irreducible 
groups $G, H\leq \glnf$ of odd order are conjugate. They isolate 
the subgroup $G^{\mathrm{Aut} G}$ of $G$ whose elements are fixed 
under $\mathrm{Aut}(G)$, and show that (with caveats) the number 
of $\GL(n,\F)$-conjugacy classes of subgroups $H$ of $\glnf$ 
isomorphic to $G$ is equal to the number of equivalence classes 
of faithful irreducible $\F$-representations of 
$G^{\mathrm{Aut} G}$.

\section{Related work}

Linear group classifications of the kind advocated by Laci are 
constantly in demand, and hence worth pursuing. The area
is effectively still wide open. We review a sample of other 
classifications with a computational aspect that relate to Laci's 
concerns. 

Dixon and B.~Mortimer~\cite{DixonMortimer} listed all primitive
permutation groups of degree less than $1000$ with insoluble 
socle. Short originally aimed to match this 
range of degrees for soluble groups. B.~Eick and 
H\"{o}fling~\cite{EickHoefling} extended the degree bound 
far beyond that in \cite{Short}. Taking Aschbacher's 
categorisation of potentially maximal subgroups of $\GL(n,q)$ 
as a starting point, they developed an algorithm to classify the
soluble irreducible subgroups of $\GL(n,p)$ for $p^n\leq 6560$.
(H\"{o}fling later augmented the list, going up to permutation
group degree $10000$.) Their algorithm involves testing irreducibility  
and subgroup conjugacy in $\GL(n,q)$. The complexity of this algorithm 
grows with $n$ and $q$. More recently, H.~J.~Coutts, M.~Quick, and 
C.~Roney-Dougal~\cite{Couttsetal} classified the insoluble irreducible 
subgroups of $\GL(n,p)$ for $p^n< 4096$.

Building on \cite{Flannery2}, and at the instigation of Laci, 
Flannery gave methods to classify the irreducible monomial 
subgroups of $\GL(4,q)$. Subsequently Flannery and 
E.~A.~O'Brien~\cite{FlanneryOBrien} designed algorithms to list 
irreducible linear groups of small degree over finite fields,
with analogous classifications over $\C$ as an ingredient. 
The input field size bounds the integer parameter strings that 
define generating sets. A key theorem in \cite{FlanneryOBrien} 
asserts that if $\F$ is any extension of $\GF(q)$, and 
$n\geq \abk 3$ or $q>\abk 3$, then a subgroup of 
$\glnf$ isomorphic to $\SL(n,q)$ is irreducible and conjugate to 
$\SL(n,q)$. The proof of this result is mainly due to Laci and uses 
\cite{Ko1}. Implementations of the algorithms of 
\cite{FlanneryOBrien} are available in {\sc Magma}.
Their efficiency depends on field arithmetic 
(as do, e.g., the algorithms of \cite{Nilmat,Glasgow}), and 
they avoid testing irreducibility or conjugacy.
The input field size is unrestricted except for a tiny number
of exceptions. This type of implementation may be compared with 
data libraries such as \cite{Short}, and with the approach 
of \cite{EickHoefling}, which requires non-trivial computation in 
$\GL(n,q)$. Techniques similar to those in \cite{FlanneryOBrien} 
could be applied at least up to degree $5$.

Classifying irreducible soluble linear groups over a finite field 
in other special degrees, such as the product of two
primes, is feasible. Suprunenko's book~\cite{SuprPermGp} is yet 
to be exploited for this purpose.
The resultant algorithms would be practical for large 
degrees and fields.

Finally, we note that progress in computational 
representation theory affords new avenues for classifying insoluble 
linear groups over finite fields and $\C$.

\section{Concluding remarks}

The second author remembers hours 
spent each week with Laci, seeing the solution of a problem 
through to its very end. Laci was extremely generous in 
sharing his expertise. The work that he inspired and
nurtured in others forms an important part of 
his legacy. We owe him a lasting debt.

\medskip

\subsubsection*{Acknowledgment}
The authors were supported by the Irish Research Council (grant `MatGpAlg')
and Science Foundation Ireland (grant 11/RFP.1/MTH/3212).
A.~S.~Detinko was further supported by a Marie Sk\l odowska-Curie Individual 
Fellowship grant H2020 MSCA-IF-2015, No.~704910 (EU Framework Programme 
for Research and Innovation). 

\bibliographystyle{line}

\bigskip

\end{document}